\documentstyle[12pt]{article}
\begin{document}
\newtheorem{theorem}{Theorem}
\newtheorem{definition}{Definition}
\newtheorem{proposition}{Proposition}
\newtheorem{lemma}{Lemma}
\begin{center}{\Large Some Conformally Flat Spin Manifolds, Dirac Operators and Automorphic Forms}\end{center}
\begin{center}{\large R. S. Krau{\ss}har,\\Department of Mathematical Analysis, Ghent University,
Galglaan 2, B-9000 Ghent, Belgium, \\and John Ryan\\Department of Mathematics, University of Arkansas, Fayetteville, \\AR 72701, USA.}\end{center}
\begin{quote}
\begin{abstract} 
In this paper we study Clifford and harmonic analysis on some conformal flat spin manifolds. In  particular we treat manifolds that can be parametrized by $U / \Gamma$ where $U$ is a simply connected subdomain of either $S^{n}$ or $R^{n}$ and $\Gamma$ is a Kleinian group acting discontinuously on $U$. Examples of such manifolds treated here include $RP^{n}$ and $S^{1}\times S^{n-1}$. Special kinds of Clifford-analytic automorphic forms associated to the different choices of $\Gamma$ are used to construct Cauchy kernels, Cauchy Integral formulas, Green's kernels and formulas together with Hardy spaces and Plemelj projection operators for $L^{p}$ spaces of hypersurfaces lying in these manifolds.
\end{abstract}
\end{quote}
\section{Introduction}
While the classical Cauchy-Riemann operator pervades much of modern analysis, particularly classical harmonic analysis, a less well known but significantly powerful differential operator arising in modern analysis is the Euclidean Dirac operator $D$. This operator has proved extremely powerful in tackling a number of problems arising in harmonic analysis and related fields. See for instance \cite{lmq,lms,m,q,q1}. The natural setting to understand the basic properties of this operator and its applications is the setting of Clifford analysis.

\ Solutions to the Dirac equation are called Clifford holomorphic functions or monogenic functions. Such functions are covariant under M\"{o}bius transformations acting over $R^{n}\cup\{\infty\}$. In fact this covariance is an automorphic invariance best described by a method using Clifford algebras and due to Ahlfors \cite{a} and Vahlen \cite{v1}. In fact all solutions to the equation $D^{k}f=0$ where $k\in {\bf N}$ exhibit a similar automorphic invariance under M\"{o}bius transformations. Given this natural method for describing conformal or M\"{o}bius transformations and the automorphic invariance of solutions to $D^{k}f=0$ under actions of the M\"{o}bius group a natural choice of generalization of Riemann surfaces from one complex variable to the present context would be conformally flat manifolds. This has previously been pointed out in \cite{ma, r85}. Conformally flat manifolds are manifolds with atlases whose transition functions are M\"{o}bius transformations.  These types of manifolds have been studied in a number of contexts independent of Clifford analysis, see for example \cite{bog,Chang,sy}. As pointed out in \cite{sy} one fruitful way of constructing conformally flat manifolds is to factor out a simply connected subdomain $U$ of either the sphere $S^{n}$ or $R^{n}$ by a Kleinian subgroup $\Gamma$ of the M\"{o}bius group where $\Gamma$ acts discontinuously on $U$. This gives rise to the conformally flat manifold $U / \Gamma$. Examples of such manifolds include for example $n$-tori, cylinders, real projective space and $S^{1}\times S^{n-1}$.

\ Our aim in this paper is to establish basic tools for developing Clifford analysis and harmonic analysis over some examples of these types of conformally flat manifolds. 

\ Monogenic, harmonic and  $k$-monogenic functions on manifolds $M$ that are parametrized by $U / \Gamma$ can be obtained from automorphic forms related to $\Gamma$ that are monogenic ($k$-monogenic) in $U$. The fundaments of a general theory of monogenic and $k$-monogenic automorphic forms associated to arithmetic hypercomplex generalizations of the modular group and its subgroups has been developed  in \cite{KraCV2,KraCMFT2}. The Eisenstein and Poincar\'e type series from  \cite{KraCV2,KraCMFT2} induce monogenic ($k$-monogenic) functions on those manifolds that are parametrized precisely by the corresponding arithmetic groups.

\ In order to develop now techniques for the treatment of boundary value problems related to  Hardy spaces on $U /\Gamma$ that arise in the context of Clifford analysis on conformally flat manifolds one needs to investigate which  automorphic forms give rise to Cauchy or Green kernel functions on $U/\Gamma$.

\ The cases of $n$-tori and cylinders which are realized by translation groups equipped with trivial bundles have already been treated in \cite{KraRyan1}.

\ Here we will deal with further important particular examples including real projective space, $RP^{n}$, and $S^{1}\times S^{n-1}$. We also deal in the cases of the tori, cylinders and real projective spaces with several examples of spinorial bundles. Furthermore, in the cases of real projective space and $S^{1}\times S^{n-1}$ we present examples of hypersurfaces for which the usual Hardy $p$-space decomposition of the $L^{p}$ space of the hypersurface breaks down. In this paper we primarily concentrate on constructing Cauchy kernels for the various conformally flat manifolds that we consider. Once one has the Cauchy kernel, the corresponding Szeg\"{o}, Poisson and Kerzman-Stein kernels and their basic properties follow along from lines layed out in \cite{KraRyan1}. For this reason we do not go into these details here.

\ Many aspects of Clifford analysis have already been developed in the context of general Riemannian manifolds, see \cite{ca,cn,mi}, including Cauchy integral formulas and Plemelj projection operators. However so far little in the way of explicit formulas have been developed in this general context. In contrast for conformally flat manifolds there is a real hope that one can introduce explicit formulas and solutions for a very wide class of settings. These include Cauchy and Green's formulas, Plemelj formulas and Kerzman-Stein kernels.

{\bf{Acknowledgement:}} The authors are grateful to David Calderbank and Michael Eastwood for many helpful discussions and e-discussions which proved to be very helpful in developing this paper.

\section{Preliminaries}

\ For particular details about Clifford algebras see for instance \cite{p} and basic Clifford analysis see for example \cite{bds}. Throughout this paper the canonical basis of the Euclidean space $R^{n}$ shall be denoted by $e_1, e_2,\ldots, e_n$. The symbol  ${\it Cl}_n$ stands for the associated real $2^{n}$-dimensional Clifford algebra in which $e_i e_j + e_j e_i = - 2 \delta_{ij}$ holds, where  $\delta_{ij}$ is the usual Kronecker symbol. Under this rule of multiplication each non-zero vector $x\in R^{n}$ has a multiplicative inverse $x^{-1}=\frac{-x}{\|x\|^{2}}$. Up to a sign this is the Kelvin inverse of the vector $x$.
We also need the anti-automorphism $\sim :Cl_{n}\rightarrow Cl_{n}:\sim e_{j_{1}}\ldots e_{j_{r}}=e_{j_{r}}\ldots e_{j_{1}}$. We shall write $\tilde{X}$ for $\sim X$.\\
The Euclidean Dirac operator is $D: = \sum\limits_{j=1}^n \frac{\partial }{\partial x_j} e_j$
and differentiable functions defined in open subsets of $R^n$ with values in ${\it Cl}_n$ that are annilhilated from the left (right) by the Dirac operator are called left (right) monogenic functions, or left (right) Clifford holomorphic functions.
The left and right fundamental solution to the $D$-operator is the euclidean Cauchy kernel function 
$G(x-y) = \frac{x-y}{\|x-y\|^n}$.
The Dirac operator factorizes the Laplacian, viz $D^2 = - \Delta$. Functions that are annihalated by the $k$-th iterate, $D^{k}$, of the Dirac operator from the left are called left $k$-monogenic functions, and right $k$-monogenic functions if they are annihilated by $D^{k}$ acting from the right. The left and right fundamental solution of $D^k$ is given by 
$$
G_k(x-y) = \left\{ \begin{array}{cc} \frac{x-y}{\|x-y\|^{n-k}}  & k \mbox{ odd integer with } k \le n-1,\\
\frac{1}{\|x-y\|^{n-k-1}}  & k \mbox{ even integer with }k \le n-1 .
\end{array} \right.
$$
These kernel functions serve as Cauchy and Green kernels in $R^{n}$ and can be used to solve special boundary value problems. See for instance \cite{m} for more details.

\ For all that follows it is crucial that the operators $D^k$ are all invariant up to a conformal or automorphic weight factor under all M\"{o}bius transformations.

\ In \cite{a,v1} and elsewhere it is shown that any M\"{o}bius transformation $\psi(x)$ over $R^{n}\cup\{\infty\}$ can be written as $y=(ax+b)(cx+d)^{-1}$ where the coefficients $a,b,c,d$ are all products of vectors from $R^{n}$, and they satisfy additionally 
$a\tilde{c}$, $c\tilde{d}$, $d\tilde{b}$ and $b\tilde{a} \in R^{n}$, and $a\tilde{d}-b\tilde{c} = \pm1$. For details see \cite{a,v1} and elsewhere.
 
\ If $f$ is a left $k$-monogenic function in the variable $y=\psi(x)=(ax+b) (cx+d)^{-1}$, then the function $J_{k}(\psi,x)f((ax+b) (cx+d)^{-1})$ is again left $k$-monogenic, now with respect to the variable $x$. $J_{1}(\psi,x)=J(\psi,x)=\frac{\widetilde{cx+d}}{\|Cx+d\|^{n}}$ and $J_{2}(\psi,x)=\frac{1}{\|cx+d\|^{n-2}}$. Precise formulas for the other $J_{k}$'s are given in \cite{r,r2} and elsewhere. 

\ We conclude this section by pointing out that as shown in \cite{r} via a Cayley transformation one can also develop Clifford analysis over the $n$-sphere, $S^{n}$. In this case the Dirac operator is $D_{s}=x(\Lambda+\frac{n}{2})$ where $x\in S^{n}$, so we are considering $S^{n}$ as embedded in the usual way in $R^{n+1}$, and $\Lambda=\Sigma_{i<j,i=1}^{n+1}e_{i}e_{j}(x_{i}\frac{\partial}{\partial x_{j}}-x_{j}\frac{\partial}{\partial x_{i}})$. In  \cite{r,v} it is shown that the Cauchy kernel for $D_{s}$ is $G_{s}(x,y)=\frac{x-y}{\|x-y\|^{n}}=\frac{x-y}{(|2-2<x,y>|)^{\frac{n}{2}}}$ where $x$ and $y\in S^{n}$. Further in \cite{lr} it is shown that the Laplacian $\triangle_{s}$ for $S^{n}$ is $D_{s}(D_{s}+x)$ and the fundamental solution to this operator is $H_{s}(x,y)=\frac{1}{n-2}\frac{1}{\|x-y\|^{n-2}}$ provided $n>2$. Again $x$ and $y\in S^{n}$. These kernels are used in \cite{lr} to produce Cauchy integral formulas and Green's formulas. See also \cite{r4,v} for the case of Cauchy integral formulas.

\section{Construction of Some Conformally Flat Spin Manifolds}

\ Conformally flat manifolds are in general $n$-dimensional manifolds that possess atlases whose transition functions are M\"obius transformation. 

\ A very systematic and constructive method to obtain conformally flat manifolds is to make use of a projection argument mentioned for instance in \cite{sy}. Suppose $U$ is a simply connected domain on either the sphere, $S^{n}$ or in $R^{n}$ and $\Gamma$ is a Kleinian group that acts discontinuously on $U$ then the factorized space $U /\Gamma$ is a conformally flat manifold.

\ This method turns out to have one great advantage.  It enables one to carry out relatively easily  Clifford and harmonic analysis techniques that have been developed previously in Euclidean spaces or on spheres and hyperbolae to a large variety of manifolds.  

\ Any M\"{o}bius transformation $y=\psi(x)$ can be either written as $(ax+b)(cx+d)^{-1}$ or as $(-ax-b)(-cx-d)^{-1}$. So in fact the automorphic invariance of Clifford holomorphic functions that we previously mentioned is correct up to a sign. So the left Clifford holomorphic function $f(y)$ is changed to $\pm J(\psi,x)f(\psi(x))$. This has an effect on constructing spinor bundles over a conformally flat manifold. Suppose that $M$ is a conformally flat manifold and $\mu_{2}\mu_{1}^{-1}:U_{1}\rightarrow U_{2}$ is a transition function arising from the atlas of $M$. So $U_{1}$ and $U_{2}$ are domains in $R^{n}$ and $\psi=\mu_{2}\mu_{1}^{-1}$ is a M\"{o}bius transformation. Let us now consider the two bundles $U_{1}\times Cl_{n}$ and $U_{2}\times Cl_{n}$. Given $u=\psi(x)\in U_{2}$ and $X\in Cl_{n}$ the pair $(u,X)\in U_{2}\times Cl_{n}$ may be identified with either $(x,J(\psi,x)X)$ or $(x,-J(\psi,x)X)$ in $U_{1}\times Cl_{n}$. If we can choose a suitable collection of signs on these local bundles which are globally compatable over $M$ then we have constructed a spinor bundle $E$ over $M$. In this case $M$ is called a conformally flat spin manifold. Note that it might be the case that there are several choices of  spinor bundles over $M$. Further it should be recalled, \cite{p}, that $Cl_{n}$ is the direct sum of several isomorphic minimal left ideals. These are often called spinor spaces. So in our construction of spinor bundles one might want to repace the Clifford algebra with one of these spinor spaces.  Following \cite{ma,r85} we may now talk about Clifford holomorphic sections. 
\begin{definition}
Given a conformally flat spin manifold $M$ with spinor bundle $E$ then a section $f:M\rightarrow E$ is called a left Clifford holomorphic section if locally $f$ reduces to a left Clifford holomorphic function. 
\end{definition}
In the previous definition the conformal weight functions $J(\psi,x)$ are used to preserve Clifford holomorphy. A similar definition can be given for right Clifford holomorphic sections over $M$. It is in this sense that conformally flat spin manifolds are natural generalizations to $n$ real dimensions of Riemann surfaces.

\ We now turn to construct some examples of conformally flat spin manifolds. Clearly $R^{n}$ is one such example. Via Cayley transformations so is $S^{n}$, \cite{r4}. Also if $\Gamma$ is a translation group acting on  $R^{n}$, then $R^{n}/\Gamma$ is a conformally flat spin manifold. In these cases we get either cylinders or $n$-tori. The links between Clifford analysis, harmonic analysis and cylinders and tori has been developed in \cite{KraRyan1}. Here are some other examples. First let $U=S^{n}$ and $\Gamma=\{\pm 1\}$ then $U / \Gamma=RP^{n}$, real projective space. Note that when $n$ is even $RP^{n}$ is no longer an orientable manifold. Now suppose $U=R^{n}\backslash\{0\}$ and $\Gamma$ is a discrete subgroup of the orthogonal group. Then $M=(R^{n}\backslash\{0\})/ \Gamma$ is a conformally flat spin manifold. Now let $U=R^{n}\backslash\{0\}$ and $\Gamma=\{2^{k}\}_{k\in{\bf{Z}}}$. Then $\Gamma$ is a discrete subgroup of the dilation group and $U / \Gamma$ is $S^{1}\times S^{n-1}$. Next let us consider a lattice $\Pi$ in $R^{n}$. Let $U=R^{n}\backslash\{0\}$ and $\Gamma=\{\psi_{m}(x)=x(mx+1)^{-1}:m\in\Pi\}$. Then let $M=U / \Gamma$.

\ In all that follows $U$ will be a universal covering space of a conformally flat manifold $M$. So there is a projection map $p:U\rightarrow M$. Further for each $x\in U$ we shall denote $p(x)$ by $x'$. Furthermore if $Q$ is a subset of $U$ then we denote $p(Q)$ by $Q'$.

\section{Some Cauchy and Green's Kernels for Some Conformally Flat Spin Manifolds}

\ We begin with real projective space, $RP^{n}$. Here $U=S^{n}$ and $\Gamma=\{\pm1\}$. We first set up a spinor bundle $E_{1}$ over $RP^{n}$ by making the identification of $(x,X)$ with $(-x,X)$ where $x\in S^{n}$ and $X\in Cl_{n}$. Now we need to change the spherical Cauchy kernel $G_{s}(x,y)$ into a kernel which is invariant  with respect to $\{\pm 1\}$ in the variable $x\in S^{n}$. We obtain $G_{s}(x,y)+G_{s}(-x,y)$. This projects to give a kernel $G_{RP,1}(x',y')$ for $RP^{n}$. Suppose now that $S$ is a suitably smooth hypersurface lying in the northern hemisphere of $S^{n}$. Suppose also that $V$ is a domain lying in the northern hemisphere and that $S$ bounds a subdomain $W$ of $V$ and $y\in W$. If $f:V\rightarrow Cl_{n}$ is a left spherical Clifford holomorphic function then 
\[f(y)=\frac{1}{\omega_{n}}\int_{S}(G_{s}(x,y)+G_{s}(-x,y))n(x)f(x)d\sigma(x)\]
where $\omega_{n}$ is the surface area of the unit sphere in $R^{n}$, and $n(x)$ is the unit outer normal vector to $S$ at $x$ lying in the tangent space of $S^{n}$ at $x$. Also $\sigma$ is the usual Lebesgue measure on $S$. As we restricted attention purely to the northern hemisphere we may now use the projection map $p:S^{n}\rightarrow RP^{n}$ to note that this projection map induces a function $f':V'\rightarrow E_{1}$. We now have
\[f'(y')=\frac{1}{\omega_{n}}\int_{S'}G_{RP,1}(x',y')dp(n(x))f'(x')d\sigma'(x'),\]
where $x'$ and $y'$ are the projections of $x$ and $y$ respectively, and $S'$ is the projection of $S$. Further our projection induces a measure $\sigma'$ on $S'$ from the measure $\sigma$ on $S$. Also $dp$ is the derivative of $p$.\\
Let us now make the situation slightly more complicated. We shall still assume that the hypersurface $S$ lies in the northern hemisphere. However, now we will assume that the domain $V$ is such that $-x\in V$ for each $x\in V$ and the spherical left Clifford holomorphic function $f$ is two fold periodic, so that $f(x)=f(-x)$. Now the projection map $p$ gives rise to a well defined domain $V'$ on $RP^{n}$ and a well defined function $f':V'\rightarrow E_{1}$ such that $f'(x')=f(\pm x)$ for $p(\pm x)=x'$. As the function $f$ is spherical left Clifford holomorphic then this construction induces a Dirac operator $D_{RP^{n}}$ on $RP^{n}$ and $D_{RP^{n}}f'(x')=0$. We shall call such functions real projective left Clifford holomorphic functions. One may similarly construct real projective right Clifford holomorphic functions. In this context one also has
\[f'(y')=\frac{1}{\omega_{n}}\int_{S'}G_{RP,1}(x',y')dp(n(x))f'(x')d\sigma'(x').\]
One way to construct a left spherical Clifford holomorphic function $f$ satisfying $f(x)=f(-x)$ is as follows. Consider a closed subset $C$ of the sphere. On this closed set one obtains a sigma algebra $\Psi$ such that for each $K\subset C$ with $K\in\Psi$ then $-K\in\Psi$. One now introduces a $Cl_{n}$ valued measure $\mu$ on $\Psi$ with the property $\mu(K)=\mu(-K)$ for each $K\in\Psi$. In this case the integral $\int_{C}(G_{s}(x,y)+G_{s}(-x,y))d\mu(y)$ defines a spherical left Clifford holomorphic function $f$  satisfying $f(x)=f(-x)$ on $S^{n}\backslash C$.\\
If now we assume that the hypersurface $S$ is such that $-S=S$ then both $y$ and $-y$ belong to the subdomain $V$ and in this case
\[\frac{1}{\omega_{n}}\int_{S'}G_{RP,1}(x',y')dp(n(x))f'(x')d\sigma(x')=2f'(y').\]
Let us now assume that $S$ is strongly Lipschitz, that $S=-S$ and that $\eta:S\rightarrow Cl_{n}$ belongs to  $L^{p}(S,Cl_{n})$ is such that $\eta(x)=\eta(-x)$. We also assume that $1<p<\infty$. Let us consider a piecewise $C^{1}$ path $y(t)\in V$ which approaches $w\in S$ nontangentially as $t$ tends to $1$. In this case 
\[\lim_{t\rightarrow 1}\frac{1}{\omega_{n}}\int_{S}(G_{s}(x,y(t))+G_{s}(-x,y(t)))n(x)\eta(x)d\sigma(x)\]
evaluates almost everywhere to 
\[\frac{1}{2}\eta(w)+P.V.\frac{1}{\omega_{n}}\int_{S}(G_{s}(x,w)+G_{s}(-x,w))n(x)\eta(x)d\sigma(x).\]
When we turn to real projective space we are now forced to consider two paths $y(t)$ and $-y(t)$ on $S^{n}$. In this case one may determine that 
\[\lim_{t\rightarrow 1}\frac{1}{\omega_{n}}\int_{S'}G_{RP,1}(x',y'(t))dp(n(x))\eta'(x')d\sigma(x')\]
evaluates to  
\[2P.V.\frac{1}{\omega_{n}}\int_{S'}G_{RP,1}(x',w')dp(n(x))\eta'(x')d\sigma(x'),\]
where $w'$ is the projection of $w$ to $RP^{n}$ and $\eta'(x')=\eta(x)$. So the usual Hardy space decomposition of $L^{p}$ spaces of hypersurfaces one sees in euclidean Clifford analysis, spherical Clifford analysis and elsewhere, see for instance \cite{lms,m}, does not always occur in case of real projective space. However it is straightforward to deduce that when $S$ lies in a hemisphere then one would get the usual Hardy space decomposition for $L^{p}(S',E_{1})=\{\eta':S'\rightarrow E_{1}$ such that $\eta'$ is $L^{p}$ integrable$\}$.

\ Besides the spinor bundle $E_{1}$ we can construct a second spinor bundle $E_{2}$ over $RP^{n}$ by identifying the pair $(x,X)$ with $(-x,-X)$ where again $x\in S^{n}$ and $X\in Cl_{n}$. In this case we need that our Cauchy kernel be antiperiodic with respect to $\Gamma=\{\pm 1\}$. So the projection map, $p$, induces a Cauchy kernel $G_{RP,2}(x',y')$ but now from the kernel $G_{s}(x,y)-G_{s}(-x,y)$.

 \ In this case a left Clifford holomorphic section $f':V'\rightarrow E_{2}$ will lift to a left Clifford holomorphic function $f:V\rightarrow Cl_{n}$ satisfying $f(x)=-f(-x)$. To construct examples of such functions recall that earlier we introduced a measure space $(C,\Psi,\mu)$ where $C$ is a closed subset of $S^{n}$. In this case the convolution $\int_{C}(G_{s}(x,y)-G_{s}(-x,y))d\mu(y)$ defines a left spherical Clifford holomorphic function $f$ satisfying $f(x)=-f(-x)$ on $S^{n}\backslash C$.

\ Suppose now that $V$ as before is a domain on $S^{n}$ and $S$ is a hypersurface in $V$ bounding a subdomain $W$ of $V$. Suppose further that $f:V\rightarrow Cl_{n}$ is a spherical left Clifford holomorphic function satisfying $f(x)=-f(-x)$. If $S$ lies entirely in one hemisphere then
\[\frac{1}{\omega_{n}}\int_{S}(G_{s}(x,y)-G_{s}(-x,y))n(x)f(x)d\sigma(x)=f(y)\] 
for each $y\in W$. Via the projection $p$ this integral formula induces the following integral formula
\[f'(y')=\frac{1}{\omega_{n}}\int_{S'}G_{RP,2}(x',y')dp(n(x))f'(x')d\sigma'(x').\]
On the other hand if $S$ is such that $S=-S$ then 
\[\frac{1}{\omega_{n}}\int_{S}(G_{s}(x,y)-G_{s}(-x,y))n(x)f(x)d\sigma(x)=0.\]
In this case the projection map $p$ gives rise to an integral over $S'$ that evaluates to zero.

\ It is now a simple exercise to determine the following. Suppose $S$ is a strongly Lipschitz hypersurface in $S^{n}$ and $S=-S$. Furthermore suppose that $\eta\in L^{p}(S)=\{\eta:S\rightarrow Cl_{n}:\|\eta\|_{p}<\infty\}$ for some $p\in (1,\infty)$. Suppose that $S$ bounds a domain $W$ and $y(t)$ is a path in $W$ with non-tangential limit $w\in S$, so $\lim_{t\rightarrow 1}y(t)=w$. Then
\[\lim_{t\rightarrow 1}\frac{1}{\omega_{n}}\int_{S'}G_{RP,2}(x',y'(t))dp(n(x))\eta'(x')d\sigma'(x')=\eta(w)\]
almost everywhere. It follows that for the bundle $E_{2}$, the hypersurface $S$ and for $1<p<\infty$ we can find for each $\eta'\in L^{p}(S',E_{2})$ a Clifford holomorphic section $f:V'\rightarrow E_{2}$ with trace, or boundary value $\eta$ almost everywhere.
 
\ It is an easy matter to determine that the Green's kernel for $E_{1}$ is obtained by applying the projection $p$ to the kernel $H_{s}(x,y)+H_{s}(-x,y)$ while for $E_{2}$ we apply the projection to $H_{s}(x,y)-H_{s}(-x,y)$. From these kernels and the Cauchy kernels on $RP^{n}$ one can readily set up a Green's formula over domains in $RP^{n}$. We leave this as a simple exercise.

\ Motivated by the construction of $E_{2}$ let us now return to the $k$ cylinder $C_{k}$ and construct a number of spinor bundles over $C_{k}$. The conformally flat manifold $C_{k}$ is obtained by factoring out $R^{n}$ by a $k$-lattice. To keep it simple we will choose the $k$-lattice to be ${\bf{Z}}^{k}={\bf{Z}}e_{1}+\ldots+{\bf{Z}}e_{k}$. In \cite{KraRyan1} the spinor bundle over $C_{k}$ is chosen to be the trivial one $C_{k}\times Cl_{n}$. However there are $k$ other spinor bundles over $C_{k}$. We shall now construct them. First let $l$ be an integer in the set $\{1,\ldots,k\}$, and consider the lattice ${\bf{Z}}^{l}={\bf{Z}}e_{1}+\ldots+{\bf{Z}}e_{l}$. there is also the lattice ${\bf{Z}}^{k-l}={\bf{Z}}e_{l+1}+\ldots+{\bf{Z}}e_{k}$. In this case ${\bf{Z}}^{k}=\{\underline{m}+\underline{n}:\underline{m}\in {\bf{Z}}^{l}$ and $\underline{n}\in {\bf{Z}}^{k-l}\}$. suppose now that $\underline{m}=m_{1}e_{1}+\ldots+m_{l}e_{l}$. Let us now make the identification $(x,X)$ with $(x+\underline{m}+\underline{n},(-1)^{m_{1}+\ldots+m_{l}}X)$ where $x\in R^{n}$ and $X\in Cl_{n}$. This identification gives rise to a spinor bundle $E^{(l)}$ over $C_{k}$.  
From the cotangent type series  
\[\cot_{q,k,0}(x) = \sum\limits_{\underline{m} \in {\bf{Z}}^{k}} G_{q} (x+\underline{m}) \quad \quad h < n-q\]
which converge normally on $R^{n}\backslash{\bf{Z}}^{k}$ we can readily obtain the kernel functions, simply viz, 
\[\cot_{q,k,0}(x,y)=\Sigma_{\underline{m}\in{\bf{Z}}^{k}}G_{q}(x-y+\underline{m})\]
as explained in our previous paper \cite{KraRyan1}. \\
Applying the projection $p_{k}:R^{n}\rightarrow C_{k}$ to these kernels induce kernels $\cot'_{q,k,0}(x',y')$ defined on $(C_{k}\times C_{k})\backslash diagonal(C_{k})$, where $diagonal(C_{k})=\{(x',x'):x'\in C_{k}\}$.

\ We can adapt these functions and kernels as follows. First for $k<n-q$ and $l\leq k$ we define the cotangent functions $\cot_{q,k,l}(x)$ to be 
\[\Sigma_{\underline{m}\in{\bf{Z}}^{l},\underline{n}\in{\bf{Z}}^{k-l}}(-1)^{m_{1}+\ldots+m_{l}}G_{q}(x+\underline{m}+\underline{n})\]
where $\underline{m}=m_{1}e_{1}+\ldots+m_{l}e_{l}$.
These are well defined functions on $R^{n}\backslash{\bf{Z}}^{k}$. From these functions we obtain the cotangent kernels
\[\cot_{q,k,l}(x,y)=\Sigma_{\underline{m}\in{\bf{Z}}^{l},\underline{n}\in{\bf{Z}}^{k-l}}(-1)^{m_{1}+\ldots+m_{l}}G_{q}(x-y+\underline{m}+\underline{n}).\]
Again using applying the projection map $p_{k}$ these kernels give rise to the kernels $\cot_{q,k,l}(x',y')$. 

\ As in the case of $E^{(0)}$ the Dirac operator $D$ on $R^{n}$ induces a Dirac operator acting on sections of the bundle $E^{(l)}$.  We shall denote this Dirac operator by $D_{l}$. Furthermore the $q$-th power $D^{q}$ of the Dirac operator $D$ induces a $q$-th order Dirac operator acting on sections of $E^{(l)}$. We denote this operator by $D_{l}^{q}$. When $q=2$ this operator is a spinorial Laplacian.

\begin{definition}
Suppose that $V'$ is a domain in $C_{k}$ and $f':V'\rightarrow E^{(l)}$ is such that $D_{l}^{(q)}f'=0$. Then $f'$ is called an $E^{(l)}$ left Clifford holomorphic section of order $q$.
\end{definition}

 \ A similar definition can be given for $E^{(l)}$ right Clifford holomorphic sections of order $q$. When $q=2$ such sections are harmonic sections for the given bundle.

\ Suppose now that $V$ is the inverse image of $V'$ under $p_{k}$. If $D_{l}^{(q)}f'=0$ then $f$ lifts to a function $f$ defined on $V$ and $D^{q}f=0$. Moreover $f(x+\underline{m}+\underline{n})=(-1)^{m_{1}+\ldots+m_{l}}f(x)$.
 
\begin{theorem}
Suppose that $f$ is as in the previous paragraph and $S$ is a surface lying in $V$ and bounding a subdomain $W$. Suppose also that for each $x\in W$ then $x+\underline{m}+\underline{n}$ not in $W$ for any $\underline{m}+\underline{n}\in {\bf{Z}}^{k}$. Then
\[f(y)=\frac{1}{\omega_{n}}\int_{S}\Sigma_{j=0}^{q-1}(-1)^{j}\cot_{j+1,k,l}(x,y)n(x)D^{j}f(x)d\sigma(x)\]
for each $y\in W$.  
\end{theorem}

\ By applying the projection map $p_{k}$ to the formula in the previous theorem one may induce an integral formula over $S'=p_{k}(S)$ to determine $f'(y')$ where $y'=p_{k}(y)$. In the case $q=1$ we obtain a Cauchy integral formula and for $q=2$ one obtains a Green's formula.

\ Adapting from \cite{KraRyan1} when $q=n-k$ the kernel $\cot_{n-k,k,l}(x,y)$ is defined by first introducing the subset $\Lambda_{r}$ of${\bf{Z}}^{r}$ where $\Lambda_{r}=\{m_{1}e_{1}+\ldots m_{r}e_{r}:M_{1},\ldots,m_{r}\in {\bf{Z}}$ and $m_{r}>0\}\cup\ldots\cup\{m_{1}:m_{1}e_{1}\in{\bf{Z}}$ and $m_{1}>0\}$. Then the series
\[G_{n-k}(x) + \Sigma_{\underline{m}\in\Lambda_{l}, \underline{n}\in\Lambda_{k-l}}(-1)^{m_{1}+\ldots+m_{l}}(G_{n-k} (x+\underline{m}+\underline{n})+G_{n-k}(x-\underline{m}-\underline{n}))\]
is defined to be $\cot_{n-k,k,l}(x,y)$.
\ Moreover again adapting from \cite{KraRyan1} we can take points $a$ and $b\in R^{n}\backslash {\bf{Z}}^{k}$ such that $a$ is not congruent to $b$ modulo ${\bf{Z}}^{k}$. Then one can define 
\[\cot_{n-k+1,k,l,a,b}(x,y)\] 
to be
\[ G_{n-k+1}(x-a) + G_{n-k+1}(x-b)+\Sigma_{\underline{m}\in\Lambda{l},\underline{n}\in\Lambda{k-l}}(-1)^{m_{1}+\ldots+m_{l}}(G_{n-k+1}(x-a+\underline{m}+\underline{n})\]
\[ +G_{n-k+1}(-a-\underline{m}-\underline{n})
+G_{n-k+1}(x-b+\underline{m}+\underline{n}) + G_{n-k+1}(-b-\underline{m}-\underline{n}).\]

\ In \cite{KraRyan1} we also introduced analogues of convolution operators of Calderon-Zygmund type acting on the $L^{p}$ spaces of certain special hypersurfaces in $C_{k}$ for $1<p<\infty$ and $1\leq k\leq n-1$. We also introduced analogues of operators of LMS type and specified Poisson, Szeg\"{o} and Bergman kernels together with Kerzman-Stein kernels. All of these readily carry over to the context considered here using the $E^{(l)}$ bundles. They carry over with only minor adaptations. For this reason we do not go into the details of constructing such kernels, but leave it as a simple exercise.

\ Turning now to the case where $U=R^{n}\backslash\{0\}$ and $\Gamma=\{2^{k}:k\in {\bf{Z}}\}$ the conformally flat spin manifold is in this case $S^{1}\times S^{n-1}$. To get a Cauchy kernel for this manifold let us start with the following series
\[C_{1}(x,y)=\Sigma_{k=0}^{\infty}G(2^{k}x-2^{k}y)\]
where $x,y\in R^{n}\backslash\{0\}$. Next put
\[C_{2}(x,y)=G(x)\Sigma_{k=-1}^{-\infty}(G(2^{-k}x^{-1}-2^{-k}y^{-1})G(y)\]
and then finally define $C(x,y)$ to be $C_{1}(x,y)+2^{2-2n}C_{2}(x,y)$. Note that $C(2x,2y)=C(x,y)$ so this kernel is periodic with respect to the Kleinian group $\{2^{k}:k\in {\bf{Z}}\}$. The Cauchy kernel $C'(x',y')$ for the projected Dirac operator $D'$ on $S^{1}\times S^{n-1}$ is then the projection of $C(x,y)$ on $(S^{1}\times S^{n-1})\times(S^{1}\times S^{n-1})\backslash diagonal(S^{1}\times S^{n-1})$.

\ The spinor bundle $E$ over $S^{1}\times S^{n-1}$ is constructed by identifying the pair $(x,X)$ with $(2^{k}x,X)$ for every $k\in{\bf{Z}}$ and with $x\in R^{n}\backslash\{0\}$ and $X\in Cl_{n}$.  For $V'$ a domain in $S^{1}\times S^{n-1}$ a section $f':V'\rightarrow E$ is called a left Clifford holomorphic section if $D'f=0$. A similar definition can be given for a right Clifford holomorphic section. Via the projection $p:R^{n}\backslash\{0\}\rightarrow S^{1}\times S^{n-1}$ the domain $V'$ lifts to an open set $V$ satisfying $2^{k}x\in V$ for each $k\in{\bf{Z}}$ and $x\in V$. Furthermore the left Clifford holomorphic section  $f'$ lifts to a left Clifford holomorphic function $f:V\rightarrow Cl_{n}$ satisfying $f(x)=f(2^{k}x)$ for each $x\in V$ and each integer $k$. Suppose that $C$ is a closed set lying in the interior of the annulus $A(0,1,3)=\{x\in R^{n}:1<\|x\|<3\}$. suppose also that $\mu$ is a measure supported on $C$. Then the convolution $\int_{C}C(x,y)d\mu(x)$ defines a left monogenic function $f$ on $V$ where $U=\{x\in R^{n}:2^{k}x\in A(0,1,3)\backslash C$ for some integer $k\}$. Furthermore a hypersurface $S'$ in $S^{1}\times S^{n-1}$ is called a strongly Lipschitz hypersurface if there is a strongly Lipschitz hypersurface $S$ lying in $R^{n}\backslash\{0\}$ and $p(S)=S'$.

\begin{theorem}
Suppose that $V'$ is a domain in $S^{1}\times S^{n-1}$. Suppose also that $S'$ is a strongly Lipschitz hypersurface in $V'$ and that $S'$ bounds a subdomain $W'$ and within $W'$ the domain $W'$ is contractable to a point. Then for each $y'\in W'$
\[f'(y')=\frac{1}{\omega_{n}}\int_{S'}C'(x',y')dp(n(x))f'(x')d\sigma'(x').\]
\end{theorem}

\ It is straightforward to calculate Plemelj projection operators for the $L^{p}$ space $L^{p}(S')=\{\theta :S'\rightarrow E: \|\theta\|_{p}<\infty\}$ and $1<p<\infty$. So for the type of hypersurface $S'$ described in the previous theorem and for $1<p<\infty$ we get the usual Hardy space decomposition $L^{p}(S')=H^{p}(S'^{+})\oplus H^{p}(S'^{-})$ where here we are assuming that the hypersurface $S'$ divides $S^{1}\times S^{n-1}$ into two complementary domains $S'^{+}$ and $S'^{-}$ and $H^{p}(S'^{\pm})$ are the Hardy $p$-spaces of left Clifford holomorphic sections on $S^{\pm}$ with $L^{p}$ non-tangential maximal functions defined on $S'$. 

\ We shall denote the subspace of $R^{n}$ spanned by $e_{1},\ldots,e_{n}$ by $R^{n-1}$. The subset $Q'=p(R^{n-1}\backslash\{0\})$ is homeomorphic to $S^{1}\times S^{n-2}$. The hypersurface $Q'$ divides $S^{1}\times S^{n-1}$ into two complementary domains $Q'^{\pm}$. It is an easy matter to see that the Poisson kernel $P'(x',y')$ for $Q'^{+}$ is induced via the projection of the real part of $2C(x,y)e_{n}$. Here $x'\in Q'$ and $y'\in Q'^{+}$. So for each $\psi'\in L^{p}(Q')$ we have that $\frac{1}{\omega_{n}}\int_{Q'}P'(x',y')\psi'(x')d\sigma'(x')$ is the solution to the Dirichlet problem for $Q'$ and $Q'^{+}$.

\ Furthermore if $T:L^{p}(R^{n-1})\rightarrow L^{p}(R^{n-1})$ is an operator of Calderon-Zygmund type with kernel $K(x-y)$, and again $1<p<\infty$. In this case it may be readily determined that the operator $T'$ defined by the projection of the series
\[L(x,y)=\Sigma_{k=0}^{\infty}K(2^{k}x-2^{k}y)+2^{2-2n}G(x)(\Sigma_{k=-1}^{-\infty}K(2^{-k}x^{-1}-2^{-k}y^{-1}))G(y)\]
is bounded over $L^{p}(Q')$ for $1<p<\infty$.

\ The Hardy $p$-space decomposition that we have just described does not occur for all reasonable hypersurfaces in $S^{1}\times S^{n-1}$. For instance suppose that $S=\{x\in R^{n}:\|x\|=2\}$ then the complement of $S'$, $=p(S)$, is connected in $S^{1}\times S^{n-1}$. Let us for the moment stick to this hypersurface. Consider for some $p\in\{1,\infty\}$ some function $\theta'\in L^{p}(S')$. Corresponding to this section is a $Cl_{n}$ valued $L^{p}$ integrable function $\theta$ defined on $S$. Now consider the integral $\frac{1}{\omega_{n}}\int_{S}C(x,y(t))n(x)f(x)d\sigma(x)$ where $y(t)$ is a path in $R^{n}\backslash\{0\}$ with non-tangential limit $w\in S$. The depending on which direction the path $y(t)$ approaches $S$ we get
\[\lim_{t\rightarrow 1}\frac{1}{\omega_{n}}\int_{S}C(x,y(t))n(x)\theta(x)d\sigma(x)=\pm\frac{1}{2}\theta(w)+\frac{1}{\omega_{n}}P.V.\int_{S}C(x,w)n(x)\theta(x)d\sigma(x)\]
almost everywhere. Under the projection $p$ this gives rise to a similar pair of formulas on $S'$ but as the complement of $S'$ is connected these formulas do not give rise to a Hardy $p$-space decomposition.
\ Let us now consider the case where $U=R^{n}\backslash\{0\}$ and where $\Gamma=\{x(\underline{m}x+1)^{-1}:\underline{m}\in{\bf{Z}}^{k}$ and $k \le n\}$ is the transversion group. We denote the conformally flat manifold $U / \Gamma$ by $M$. A spinor bundle over $M$ is constructed by identifying the pair $(u,X)$ with $(x,\frac{\widetilde{\underline{m}x+1}}{\|\underline{m}x+1\|^{n}}X)$ for each $\underline{m}\in{\bf{Z}}^{k}$. Here $u=x(\underline{m}x+1)^{-1}$ and as before $x\in R^{n}\backslash\{0\}$ and $X\in Cl_{n}$. Again we have a projection $p:R^{n}\backslash\{0\}\rightarrow M$ and this projection induces from the Dirac operator $D$ a new Dirac operator $D'$ acting on sections taking values in $E$. Given a domain $U'\subset M$ a section $f':U'\rightarrow E$ is called a left Clifford holomorphic section if $D'F'=0$. A similar definition can be given for right Clifford holomorphic sections. Let us now construct a Cauchy kernel for such sections. For the trivial spinor bundle in the case $k < n-1$ the Cauchy kernel is simply given by the series  
\[G(x,y):= G(x) ( \Sigma_{\underline{m} \in {\bf{Z}}^{k}} G(x^{-1}-y^{-1}+\underline{m})) G(y).\]
For the cases $k\ge n-1$ this series does not converge. In these cases take for $k=n-1$ 
$$
G(x,y):= G(x) \Big( G(x^{-1} - y^{-1}) +  \Sigma_{\underline{m} \in {\bf{Z}}^{k}} G(x^{-1}-y^{-1}+\underline{m}) - G(\underline{m})\Big) G(y),
$$
and in the case $k=n$ 
\begin{eqnarray*}
G(x,y) &:=& G(x) \Big( G(x^{-1} - y^{-1} -a) - G(x^{-1} - y^{-1} -b) + \\
& & \sum\limits_{\underline{m} \in {\bf Z} \backslash\{0\}} \Big[ G(x^{-1} - y^{-1} - a + \underline{m}) - G(x^{-1} - y^{-1} - b + \underline{m}) \\
& & \quad \quad \quad -G(\underline{m}-a) + G(\underline{m}-b)\Big] \Big) G(y).
\end{eqnarray*}
where  $a,b \in R^n\backslash {\bf Z}^n$ with $a \not\equiv b$ mod ${\bf Z}^k$. 
On placing $u=x(\underline{m}x+1)^{-1}$ it may be observed that $\tilde{u}=u$. So $u=(\widetilde{\underline{m}x+1})^{-1}\tilde{x}=(\widetilde{\underline{m}x+1})^{-1}x$. Using this observation it may be determined that 
\[\frac{\widetilde{\underline{m}x+1}}{\|\underline{m}x+1\|^{n}}G(u,v)\frac{\underline{m}y+1}
{\|\underline{m}y+1\|^{n}}=G(x,y)\]
where $v=y(\underline{m}y+1)^{-1}$. 
This last formula shows that the kernel $G(x,y)$ is an automorphic form for the transversion group. In contrast to the automorphic forms that are constructed in \cite{KraCMFT1} for this group, which have a one-sided automorphic form invariance, this kernel $G(x,y)$ has a two sided automorphic form invariance similar to that described in \cite{KraCMFT2}.  Again it may be observed that the projection map $p:R^{n}\backslash\{0\}\rightarrow M$ induces a kernel $G'(x',y')$ on $(M\times M)\backslash diagonal(M)$. From this one obtains a Cauchy integral formula similar to the one established in Theorem 2.  Notice that the series $G(x,y)$ arise from the cotangent type functions by  Kelvin transformation. Indeed, the transversion group is conjugated to the translation group.

\ Similarly, to the case dealing with the $k$-cylinders, there are $k$ other spinor bundles over $M$, namely by identifying the pair $(x,X)$ with the pair $(x((\underline{m} + \underline{n})x+1)^{-1}, (-1)^{m_1+\cdots+m_l} \frac{\widetilde{(\underline{m}+ \underline{n})x+1}}{\|(\underline{m}+\underline{n})x+1\|^{n}}X)$ for each $\underline{m}\in{\bf{Z}}^{l}$ and $\underline{n} \in  {\bf Z}^{k-l}$ where we assume $1 \le l < k$. For $k<n-1$ the associated kernel function is thus given by 
\[G(x,y):= G(x) ( \Sigma_{\underline{m} \in {\bf{Z}}^{l},\underline{n} \in {\bf Z}^{k-l}} (-1)^{m_1+\cdots+m_l} G(x^{-1}-y^{-1}+\underline{m})) G(y).\]
Similarly, as presented in the cases $k=n-1$ and $k=n$ previously for the $n-1$-cylinder and the $n$-torus, one gets an analogous representation for the Cauchy kernel within the context of the transversion group for $k=n-1$ and $k=n$.

\ We conclude by suggesting that a couple further examples can be constructed explicitely simply by taking combinations of those that we have treated so far. To be more precise, we can  form semi-direct products of the groups considered earlier and treat manifolds of the type $U/\Gamma_1 \times \Gamma_2$, where $U$ is the discontinuity domain of the group $\Gamma_1 \times \Gamma_2$. To leave it simple we restrict to  the concrete treatment of the particular example factoring $R^n  \backslash\{0\}$ by the semi-direct product of a translation group and the group $1,-1$. Indeed, the whole domain $R^n \backslash\{0\}$ is fixpoint free, so that we get actually a one to one projection down to the manifold. Combining the results that we obtained before, we get $2k$ Cauchy kernels since we have now $2k$ different spinor bundles involved.  More generally, the associated Green kernels are thus given by $G(x,y) = \cot_{q,k,l}(x,y) \pm \cot_{q,k,l}(-x,y)$ where $l \in \{1,\ldots,k\}$. In absolutely the same way one can treat other examples of manifolds that can be constructed from semi-direct products of groups treated earlier, for instance factoring $R^n\backslash\{0\}$ by the semi-direct product of the transversion group with $1,-1$, etc.. The representation formulas for the associated kernels can again be readily  deduced by combining the representation formulas described here earlier.  This method opens the door to treat a wide class of conformally flat manifolds.  


\begin{thebibliography}{999}
\bibitem{a} L. V. Ahlfors, {\it{M\"{o}bius transformations in $R^{n}$ expressed through $2\times 2$ matrices of Clifford numbers}}, Complex Variables, 5, 1986, 215-224.
\bibitem{bog} T. Branson, G. Olafsson and P. Gilkey, {\it{Invariants of conformally flat manifolds}}, transactions of the AMS, 347, 1995, 939-954.
\bibitem{bds} F. Brackx, R. Delanghe and F. Sommen, {\it{Clifford Analysis}}, Pitman 76, London, 1982.
\bibitem{ca} D. Calderbank, {\it{Dirac operators and Clifford analysis on manifolds with boundary}}, Max Planck Institute for Mathematics, Bonn, preprint number 96-131, 1996.
\bibitem{Chang} A. Chang, J. Qing and P. Yang, {\it{Compactification of a class of conformally flat 4-manifold}} Invent. Math., 142, 2000, 65-93.
\bibitem{cn} J. Cnops, {\it{An Introduction to Dirac Operators on Manifolds}}, Progress in Mathematical Physics, Birkh\"{a}user, Boston, 2002.
\bibitem{KraCV2} R. S. Krau{\ss}har, {\it{Automorphic forms in Clifford analysis}}, Complex Variables, 47
No.5, 2002, 417-440.
\bibitem{KraCMFT1} R.S. Krau{\ss}har, Eisenstein Series in Complexified Clifford Analysis, to appear in {\it Computational Methods and Function Theory}.
\bibitem{KraCMFT2} R.S. Krau{\ss}har, Monogenic Modular Forms in Two and Several Real and Complex Vector Variables, to appear in {\it Computational Methods and Function Theory}.
\bibitem{KraRyan1} R. S. Krau{\ss}har and John Ryan, {\it Clifford and Harmonic Analysis on Cylinders and Tori}, to appear.
\bibitem{lmq} C. Li, A. McIntosh and T. Qian, {\it{Clifford algebras, Fourier transforms, and singular convolution operators on Lipschitz surfaces}}, Revista Mathematica Iberoamericana, 10, 1994, 665-721.
\bibitem{lms} C. Li, A. McIntosh and S. Semmes, {\it{Convolution singular integrals on Lipschitz surfaces}}, Journal of the American Mathematical Society, 5, 1992, 455-481.
\bibitem{lr} H. Liu and J. Ryan, {\it{Clifford analysis techniques for spherical pde}} to appear in Journal of Fourier Analysis and its Applications.
\bibitem{ma} M. Markel, {\it{Regular functions over conformal quaternionic manifolds}}, Commentationes Mathematicae Universitatis Carolinae, 22, 1981, 579-583.
\bibitem{m} A. McIntosh, {\it{Clifford algebras, Fourier theory, singular integrals, and harmonic functions on Lipschitz domains}}, Clifford Algebras in Analysis and Related Topics, edited by J. Ryan, CRC Press, Boca Raton, 1996, 33-87.
\bibitem{mi} M. Mitrea, {\it{Generalized Dirac operators on nonsmooth manifolds and Maxwell's equations}}, Journal of Fourier Analysis and its Applications, 7, 2001, 207-256.
\bibitem{p} I. Porteous, {\it{Clifford Algebras and Classical Groups}}, Cambridge University Press, Cambridge, 1995.
\bibitem{q} T. Qian, {\it{Singular integrals with monogenic kernels on the $m$-torus and their Lipschitz perturbations}}, Clifford Algebras in Analysis and Related Topics, edited by J. Ryan, CRC Press, Boca Raton, 1996, 157-171.
\bibitem{q1} T. Qian, {\it {Singular integrals on star-shaped Lipschitz surfaces in the quaternionic space}}, Math. Ann. 310, No.4 1998, 601-630.
\bibitem{r85} J. Ryan, {\it{Conformal Clifford manifolds arising in Clifford analysis}}, Proc. R. Ir. Acad., Sect. A 85  1985, 1-23.
\bibitem{r} J. Ryan, {\it{Iterated Dirac operators in $C^{n}$}}, Zeitschrift f\"{u}r Analysis und ihre Anwendungen, 9, 1990, 385-401.
\bibitem{r2} J. Ryan, {\it{Conformally covariant operators in Clifford analysis}}, Zeitschrift f\"{u}r Analysis und ihre Anwendungen, 14, 1995, 677-704.
\bibitem{r4} J. Ryan, {\it{Dirac operators on spheres and hyperbolae}}, Bolletin de la Sociedad Matematica a Mexicana, 3, 1996, 255-270.
\bibitem{sy} R. Schoen and S.-T. Yau, {\it{Conformally flat manifolds, Kleinian groups and scalar curvature}}, Inventiones Mathematica, 92, 1988, 47-71.
\bibitem{v1} K. Th. Vahlen, {\it{\"{U}ber Bewegungen und Complexe Zahlen}}, Math. Ann., 55, 1902, 585-593.
\bibitem{v} P. Van Lancker, {\it Clifford analysis on the sphere}, Clifford Algebras and their Applications in Mathematical Physics, edited by V. Dietrich et al, Kluwer, Dordrecht, 1998, 201-215.
\end{thebibliography}
\end{document}